\documentclass[letterpaper, 10 pt, conference]{Styles/ieeeconf}  

\IEEEoverridecommandlockouts                              
\usepackage{multicol}
\usepackage{graphicx}        
\usepackage{amsmath,amssymb}
\usepackage{caption}
\usepackage[labelformat=simple]{subcaption}
\usepackage{Styles/romoco_ames_sty}
\usepackage{algorithmic}
\usepackage{algorithm}
\usepackage{xspace}
\usepackage[shortcuts]{extdash}

\DeclareCaptionType{copyrightbox}

\makeatletter
\providecommand*{\toclevel@algorithm}{0}

\renewcommand\p@subfigure{\thefigure}
\makeatother

\newcommand{\Y}{\mathrm{Y}_{C}}
\newcommand{\U}{\mathrm{U}_{C}}
\newcommand{\Uhat}{\hat{\mathrm{U}}_{C}}
\newcommand{\Ndof}{b}

\title{\LARGE \bf
System Identification and Control of Valkyrie through SVA--Based Regressor Computation
}
\author{Shishir Kolathaya$^{1}$, Benjamin J. Morris$^{1}$, Ryan W. Sinnet$^{1}$ and Aaron D. Ames$^{2}$
\thanks{*This research is supported by NASA grant NNX11AN06H, NSF grants CNS-0953823 and CNS-1136104, and NHARP award 00512-0184-2009.}
\thanks{$^{1}$ R. Sinnet and B. Morris are with the department of Mechanical Engineering,
        Texas A \& M University, 3128 TAMU, College Station, USA
        {\tt\small {shishirny,rsinnet,bmorris}@tamu.edu}}%
\thanks{$^{2}$ S. Kolathaya and Prof. A. D. Ames are with department of Mechanical Engineering,
        Georgia Institute of Technology, 85 5th St, Atlanta, USA
        {\tt\small {shishirny,ames}@gatech.edu}}%
}

\begin{document}

\maketitle
\thispagestyle{empty}
\pagestyle{empty}

\begin{abstract}

This paper demonstrates simultaneous identification and control of the humanoid robot, Valkyrie, utilizing Spatial Vector Algebra (SVA).
In particular, the inertia, Coriolis-centrifugal and gravity terms for the dynamics of a robot are computed using spatial inertia tensors. With the assumption that the link lengths or the distance between the joint axes are accurately known, it will be shown that inertial properties of a robot can be directly evaluated from the inertia tensor. An algorithm is proposed to evaluate the regressor, yielding a run time of $O(n^2)$. The efficiency of this algorithm yields a means for online system identification via the SVA--based regressor and, as a byproduct, a method for accurate model-based control. Experimental validation of the proposed method is provided through its implementation in three case studies: offline identification of a double pendulum and a $4$-DOF robotic leg, and online identification and control of a $4$-DOF robotic arm.

\end{abstract}

\section{Introduction}

The field of system identification has been a subject of much attention in the field of robotic systems, since the 1980's \cite{kozlowski1998modelling,bukkems2003online,spong2006robot,dhparam}. One of the most important reasons is that achieving good tracking performance in robotic systems, specifically exponential convergence, without knowing the complete model of the robot has not been shown. Asymptotic convergence in tracking without knowing the model parameters has been shown in \cite{slotine1987adaptive,ghorbel1989adaptive} by using a special property of the Lagrangian dynamics of the robot, that is, the linearity of the parameters in the dynamics. Exponential convergence specific to a particular task using machine learning was also shown in \cite{horowitz1991exponential} by using the concept of persistence of excitation \cite{moore1990functional}, which requires running a series of trials for the controller to learn.

As a deviation from the methods shown above, it can be argued that identification of physical systems is required to realize good tracking performance. 
 Currently, many system identification procedures have been implemented mainly by either measuring the parameters of the robot part by part \cite{khosla1985parameter}, or dynamically by using the acceleration, velocity and angles of the robot \cite{dhparam,pedrocchi2013robot}. \cite{atkeson1986estimation} showed that the parameters that are identified as linear combinations can be consistently set to zero to determine the set of identifiable parameters. Specifically, since the parameters are affine, the equation of motion can be expressed as a matrix (regressor) multiplied by a vector of unknown parameters (base inertial parameters). \cite{yang1986simplification} used a novel method of selectively designing a robot such that the resulting inertia distribution linearizes the manipulator dynamics. In other words, the inertia parameters are made affine in the equations of motion of the rigid body robot. \cite{sheu1991identifying} used a similar method, where the affineness of the parameters in the dynamics is leveraged to compute the regressor. 

\begin{figure}[t!]
\centering
 \begin{subfigure}[b]{0.48\columnwidth}
\includegraphics[height=1.5\columnwidth]{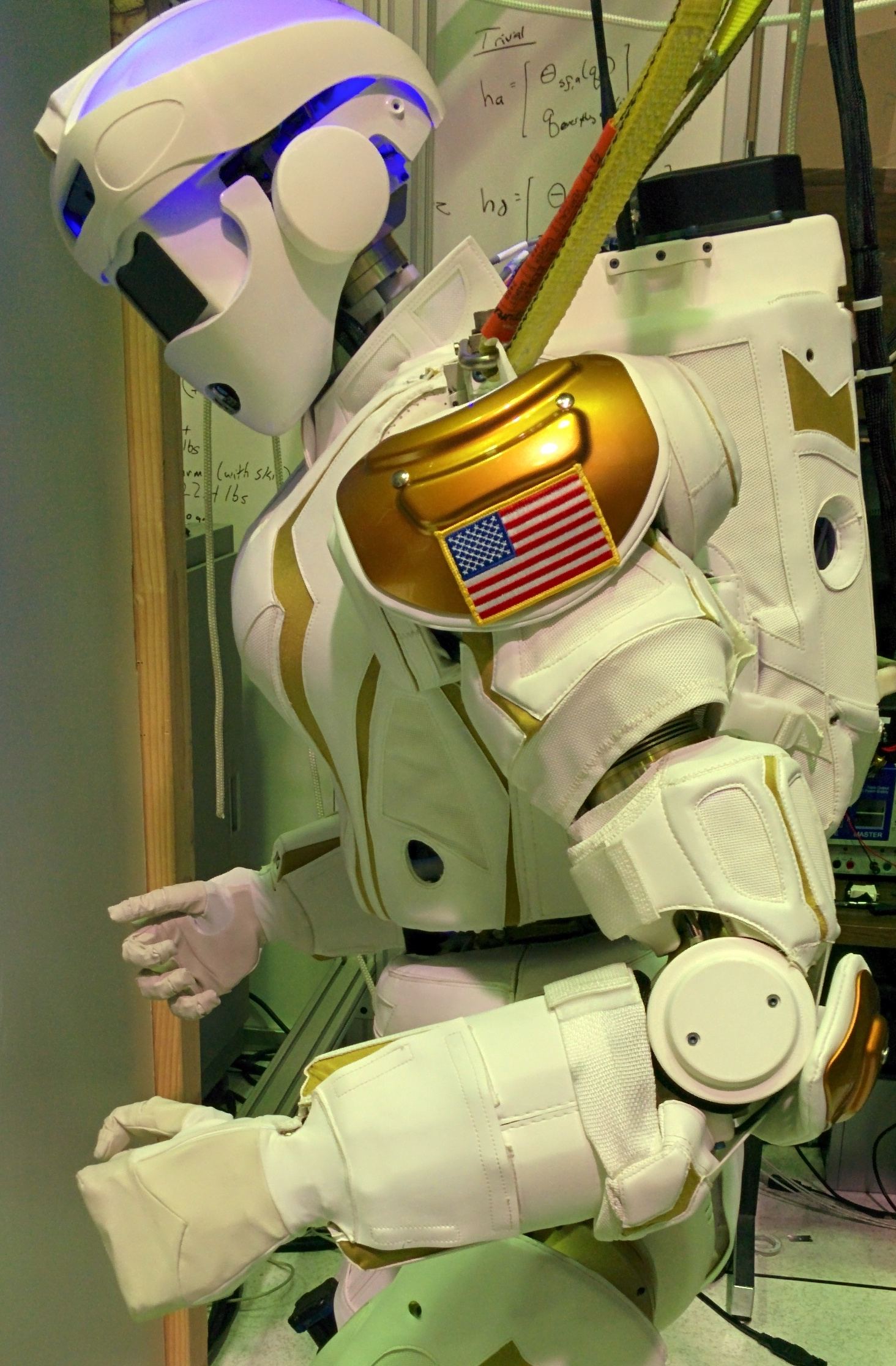}  
 \end{subfigure}
\begin{subfigure}[b]{0.48\columnwidth}
\includegraphics[height=1.5\columnwidth]{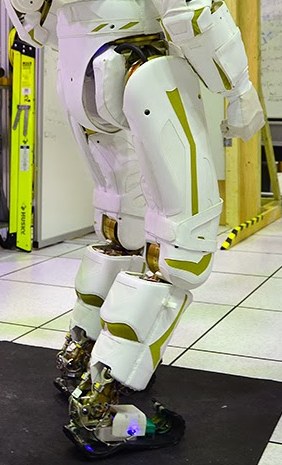}    
 \end{subfigure}
\caption{Figures showing the 4-DOF robotic arm on the left and the 4-DOF leg on the right of the Valkyrie robot on which identification was
conducted.}
\label{fig:valk}
\end{figure}
Computing the regressor is primarily done by solving for the dynamics for an $n$-DOF, $b$-body robot and collecting the unknown parameters into a vector. \cite{lu1993regressor} approached this problem using an energy based approach by using the Lagrangian formulation of robot dynamics as a starting point. \cite{lu1993regressor} also showed a second approach where the Newton-Euler recursion method is reformulated using vector analysis--type techniques. \cite{atkeson1986estimation} also used Newton-Euler equations in which the the acceleration data are obtained through a least squares estimation and the applied torques and forces are substituted to evaluate the regressor. 

This paper uses the method adopted from \cite{niemeyer1991performance} to evaluate the regressor, i.e., use Spatial Vector Algebra (SVA) to compute the regressors, an elegant way of representing the Newton-Euler equations.
If the distance between the axes are known prior to the experiment, the parameters to be identified are effectively the contents of the {\em spatial inertia tensor}. These tensors are computed by shifting the inertial elements to the joint axes. The contents of the tensor are also similar to the D-H parameters of each link in \cite{dhparam}. Given a n-DOF robot, this technique specifically lists the parameters to be identified directly from the spatial inertia tensor. Contents of this tensor are not the minimum representation, and therefore will not be unique. But, it will be shown that these non-unique parameters obtained are sufficient for realizing the model based controller, computed torque, on the robot. 
This method is demonstrated on a double pendulum as well as the leg and arm of the Valkyrie robot (\figref{fig:valk}).
Online identification is done on the pendulum and the leg, and online model based control combined with identification is done on the arm.

%
%
%

We start with a brief introduction to spatial vectors in \secref{sec:sva} which is extracted from \cite{featherstone:RBDL}. Representation of kinetic energy, spatial momentum, forces and rigid body transformations in terms of spatial vectors are also explained. \secref{sec:lag} shows how to use Spatial Vector Algebra (SVA) to extract the base inertial elements and the regressor of the robot conveniently. The derivation of this regressor and the base inertial elements are explained in detail in \secref{sec:reg}. The resulting algorithm to compute the regressor is explained in the same section and applications to control are considered in \secref{sec:estctrl}. This is finally implemented and parameters are identified for three models in \secref{sec:casestudy}.



\section{Spatial Vector Algebra for a \\Rigid Body}
\label{sec:sva}
This section will introduce the concept of spatial vectors and Spatial Vector Algebra. This is primarily derived from \cite{featherstone:RBDL} and many of the equations in this section are variants of the equations found in the same.

Rigid body motions and forces are normally described in two separate entities: 3D linear vectors and 3D angular vectors. Computing linear and rotational dynamics separately has been the normal practice in describing the equations of motion of bodies. But, if the two 3D vectors are combined together to form a 6D vector, a new vector space can be described for such systems. This vector space, of course has different rules and regulations when performing the standard mathematical operations. The 6D vectors are formed formally by what are called the Pl\"ucker coordinates.

\begin{figure}[b]
\centering
 \includegraphics[width=0.35\columnwidth]{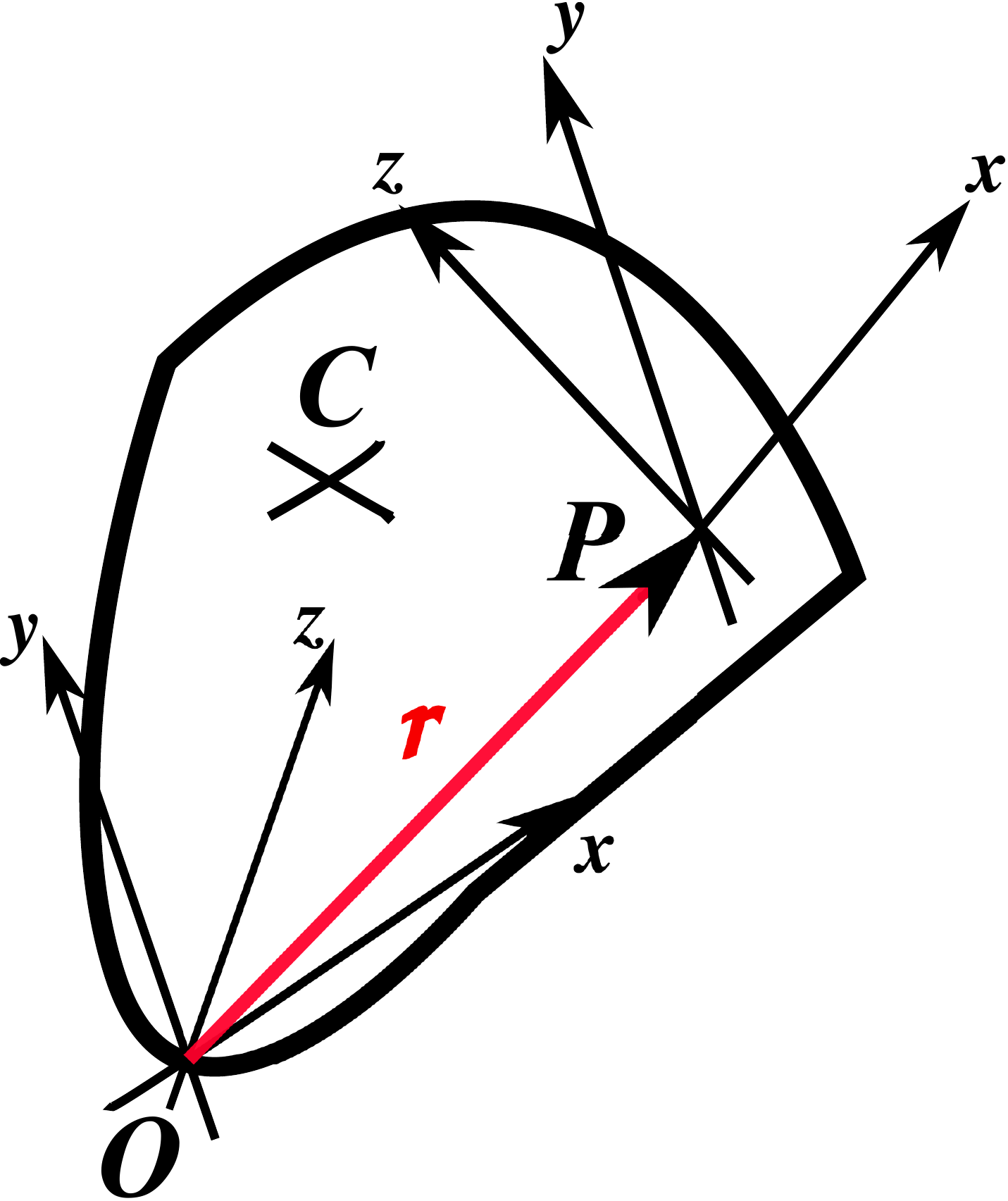}
 \caption{Figure showing the rigid body with the origin $O$, an arbitrary point $P$ and the center of mass located at $C$.}
 \label{fig:sva}
\end{figure}

A rigid body with origin located at point $O$, linear velocity $v$, and angular velocity $\omega$, about an axis passing through $O$ is shown in \figref{fig:sva}. The spatial velocity of the rigid body can be represented as:
\begin{eqnarray}
\label{eq:spatialvelocity}
\hat{v}_O = \left [ \begin{array}{c}
						\omega_x  ,
						\omega_y  ,
						\omega_z  ,
						v_x	  ,
						v_y 	  ,
						v_z		\end{array} \right ]^T
\end{eqnarray}
Similarly, the spatial force which consists of the linear force acting at point $O$ and the moment about the axis passing through $O$ can be represented as:
\begin{eqnarray}
\label{eq:spatialforce}
\hat{f}_O = \left [ \begin{array}{c}
						\tau_x  ,
						\tau_y  ,
						\tau_z  ,
						f_x		,
						f_y 		,
						f_z			\end{array} \right ]^T.
\end{eqnarray}
The order of angular and translational vectors considered is not important. Spatial vectors can also be considered with translational vector considered first and followed by the rotational vector.

It is important to note that the spatial vectors and forces are independent of the origin considered and depend solely on the bases chosen. The reason behind using spatial vectors is that both rotations and translations can be represented in one vector. In addition, the properties of these spatial vectors are different and they have a different algebra. More details about the algebra of spatial vectors can be found in \cite{featherstone:RBDL}.

\newsec{Coordinate Transforms.} It is important to have coordinate transformations in order to realize rotations and translations of coordinate frames of robotic systems using spatial vectors.

\newsec{Translation.} Spatial forces and velocities have the following forms of translation from point $O$ to an arbitrary point $P$.
\begin{eqnarray}
\label{eq:translationvelocity}
\hat{v}_P = \left [ \begin{array}{cc}
					\mathrm{1} & \mathrm{0} \\
					-r \times  & \mathrm{1} \end{array} \right ] \hat{v}_O, \:\:\:
\hat{f}_P = \left [ \begin{array}{cc}
					\mathrm{1} & -r\times \\
					\mathrm{0} & \mathrm{1} \end{array} \right ] \hat{f}_O
\end{eqnarray}
These can be derived based on the fact that the linear velocity $v_P = v_O + \omega \times r$ and torque $\tau_P= \tau_O + f \times r$. Here $r$ is the position vector directed from $O$ to $P$, i.e., $\vec{OP}$ (see \figref{fig:sva}). $r\times$ denotes the matrix equivalent of the cross product.

\newsec{Rotation.} Spatial forces and velocities have the following forms of rotation about a point $O$:
\begin{eqnarray}
\label{eq:rotation}
\hat{v}_P = \left [ \begin{array}{cc}
					\nu & \mathrm{0} \\
					\mathrm{0} & \nu \end{array} \right ] \hat{v}_O, \:\:\: \hat{f}_P = \left [ \begin{array}{cc}
					\nu & \mathrm{0} \\
					\mathrm{0} & \nu \end{array} \right ] \hat{f}_O
\end{eqnarray}
where $\nu$ indicates the rotation about the axes $x,y$ or $z$ or two of them or even all of them at once. If both rotations and translations are involved then the following relationship is obtained:
\begin{eqnarray}
\label{eq:transrotationcombined}
\hat{v}_P = \left [ \begin{array}{cc}
					\nu & \mathrm{0} \\
					-\nu r \times  & \nu \end{array} \right ] \hat{v}_O, \:\:\: \hat{f}_P = \left [ \begin{array}{cc}
					\nu & -\nu r\times \\
					\mathrm{0} & \nu \end{array} \right ] \hat{f}_O.
\end{eqnarray}
Note that the translation operation from $O$ to $P$ is done first, and the the rotation about point $P$ is carried out. If it is required that the rotation be done first, then the rotation matrix from \eqref{eq:rotation} is multiplied with $v_O$, and then the translation is carried out. In doing this, it is important to remember that the vector $r$ also gets rotated by $\nu$. The resulting coordinate transformation will look like:
\begin{eqnarray}
\label{eq:rottranscomb}
\hat{v}_P = \left [ \begin{array}{cc}
					\nu & \mathrm{0} \\
					-{(\nu r)\times} \:\: \nu & \nu \end{array} \right ] \hat{v}_O, \nonumber \\ \hat{f}_P = \left [ \begin{array}{cc}
					\nu & -{(\nu r)\times} \:\: \nu \\
					\mathrm{0} & \nu \end{array} \right ] \hat{f}_O,
\end{eqnarray}
where $(\nu r)$ is the vector $r$ rotated by the matrix $\nu$. Note that $(\nu r)\times = \nu {r\times} \nu^{-1}$. When this is substituted in \eqref{eq:rottranscomb} we effectively get \eqref{eq:transrotationcombined}. This result will be useful in reconstructing the spatial inertia tensors in order to conveniently evaluate the regressor algorithm.

\newsec{Momentum of a Rigid Body.} If the body has mass $m$, rotational inertia $\bar{I}_C$ about its center of mass, the following spatial momentum is described:
\begin{eqnarray}
\label{eq:momentum}
\hat{h}_C = \left [ \begin{array}{c} 
						\bar{I}_C \omega \\ m v_C \end{array} \right ] = \underbrace{\left [ \begin{array}{cc}
																		\bar{I}_C & \mathrm{0} \\
																		\mathrm{0} & m \mathrm{1} \end{array} \right ]}_{I_C} \hat{v}_C,
\end{eqnarray}
which is the product of the spatial inertia $I_C$ and the spatial velocity $\hat{v}_C$. The momentum defined in \eqref{eq:momentum} was w.r.t. the center of mass. To compute the momentum about an arbitrary point $O$, we have to do the transformation. So the transformation from point $O$ to the center of mass $C$ is given by:
\begin{eqnarray}
\label{eq:momentumtransformation}
\hat{h}_C = \left [ \begin{array}{cc} \nu & 0 \\ 0 & \nu \end{array} \right ] \left [ \begin{array}{cc}
\mathrm{1} & - r \times \\
\mathrm{0} & \mathrm{1} \end{array} \right ] \hat{h}_O, 
\end{eqnarray}
which is obtained since $m v_O = m v_C +m {r\times} \omega_C$. $v_O$, $\omega_O$ represent the spatial vectors at point $O$, and $v_C,\omega_C$ represent the spatial vectors at point $C$. $r$ is new position vector from point $O$ to the center of mass $C$ (instead of $P$). Expressing $\hat{h}_O$ in terms of $\hat{h}_C$, the transformation matrices get inverted:
\begin{eqnarray}
\label{eq:momentumtransformation2}
\hat{h}_O =  \left [ \begin{array}{cc}
\mathrm{1} & r \times \\
\mathrm{0} & \mathrm{1} \end{array} \right ] \left [ \begin{array}{cc} \nu & 0 \\ 0 & \nu \end{array} \right ]^{-1} \hat{h}_C, 
\end{eqnarray}
and using \eqref{eq:momentum} in \eqref{eq:momentumtransformation} and substituting for $\hat{v}_C$, we have:
\begin{eqnarray}
\hat{h}_O &=& \left [ \begin{array}{cc}
\nu^{-1} & {r \times} \: \nu^{-1}\\
\mathrm{0} & \nu^{-1} \end{array} \right ] 	I_C \hat{v}_C \\
 &=&\left [ \begin{array}{cc}
\nu^{-1} & {r \times} \: \nu^{-1}\\
\mathrm{0} & \nu^{-1} \end{array} \right ] 	I_C \left [ \begin{array}{cc}
					\nu & \mathrm{0} \\
					-\nu r \times  & \nu \end{array} \right ] \hat{v}_O. \nonumber
\end{eqnarray}
If the center of mass of the rigid body is $c$, then at zero rotation angle, let $r=c$. In other words, let $r$ be the position of the center of mass such that rotation of the coordinate frame results in negative rotation of $r$. In other words $r = \nu^{-1} c$. Applying the trick used in \eqref{eq:rottranscomb} and substituting for $r\times = (\nu^{-1}c)\times = \nu^{-1} c\times \nu$, we have the following result:
\begin{eqnarray}
\hat{h}_O =  \left [ \begin{array}{cc} \nu & 0 \\ 0 & \nu \end{array} \right ]^{-1}  I_O \left [ \begin{array}{cc} \nu & 0 \\ 0 & \nu \end{array} \right ] \hat{v}_O, 
\end{eqnarray}
where $I_O$:
\begin{eqnarray}
I_O= \left [ \begin{array}{cc}
\mathrm{1} & c \times \\
\mathrm{0} & \mathrm{1} \end{array} \right ]  I_C \left [ \begin{array}{cc}
\mathrm{1} &  \mathrm{0} \\
c \times^T & \mathrm{1} \end{array} \right ], 
\end{eqnarray}
forms the spatial inertia tensor, which is purely a function of the parameters of the robot and independent of the orientation of the coordinate frame considered. This equation will be used for a general n-DOF b-body robot where the spatial inertia tensor is effectively utilized to compute the unknown parameters. Simplifying the spatial inertia tensor results in:
\begin{eqnarray}
\label{eq:inertiatensor}
I_0=\begin{bmatrix}
		\bar{I}_C+ m \:\: {c\times} \:\: {c\times}^T & m \:\:c\times \\
		m\:\:{c\times}^T  & m\mathrm{1} \end{bmatrix}.
\end{eqnarray}

Having the expression for momentum, the kinetic energy can now be computed as:
\begin{eqnarray}
 \label{eq:kineticenergyrigidbody}
 T &=& \frac{1}{2}\hat{h}^T \hat{v} = \frac{1}{2}\hat{h}_O^T \hat{v}_O \nonumber \\
 &=& \frac{1}{2}\hat{v}_O^T \left [ \begin{array}{cc} \nu & 0 \\ 0 & \nu \end{array} \right ]^{T} I_O \left [ \begin{array}{cc} \nu & 0 \\ 0 & \nu \end{array} \right ]  \hat{v}_O.
\end{eqnarray}

\section{Lagrangian Dynamics for an \\\lowercase{n}-DOF, \lowercase{b}-Body Robot}
\label{sec:lag}

Equation of motion of an n-DOF manipulator is explained in detail in this section.

A robot can be modeled as an $n$-link manipulator. Given the configuration space $\mathbb{Q} \subset \R^n$, with the coordinates $\q \in \mathbb{Q}$, and the velocities $\dq\in T_\q \mathbb{Q}$, the Lagrangian of the $n$ degree of freedom robot can be defined as:
\begin{align}
\label{eq:Lag}
\LagB(\q, \dq) = \frac{1}{2} \dq^T D(\q) \dq - V(\q),
\end{align}
where $D(\q)\in \R^{n\times n}$ is the mass matrix of the robot, $V(q)\in \R^{n}$ is the potential energy of the robot. Specifically. The equations of motion of the $n$-link robot can be derived as:
\begin{eqnarray}
 \label{eq:eom}
 D(\q)\ddq +C(\q,\dq)\dq+G(\q) &=B\mathrm{U},
\end{eqnarray}
where $C(\q,\dq)\in \R^{n \times n}$ is the matrix of Coriolis and centrifugal forces and $G(\q)\in \R^{n}$ is the gravity matrix, $\mathrm{U} \in \R^m$ is the torque input with $m$ being the number of actuators, and $ B\in \R^{n\times m}$ is the mapping from actuator torques to joint torques, often the identity map.

The first step is to determine the inertia $D(\q)$, Coriolis-centrifugal $C(\q,\dq)$ and gravity $G(\q)$ matrices of the robot via spatial vectors. We define the body spatial velocity for each link or body $i$ of the manipulator about the joint axis $O_i$:
\begin{eqnarray}
 \label{eq:bodyvelocity}
 v_{O_i}= J_i(\q) \dq,
\end{eqnarray}
$J_i(\q)$ is the body Jacobian of link $i$ for the $i^{th}$ joint axis. Declare the spatial inertia tensor about the joint axis $O_i$ for the $i^{th}$ link or body as $I_{O_i}$ which is obtained from \eqref{eq:inertiatensor}:
\begin{eqnarray}
\label{eq:inertiatensorspecific}
{I}_{O_i}=\begin{bmatrix}
		\bar{I}_{C_i}+ m_i \:\: {c_i\times} \:\: {c_i\times}^T & m_i \:\: c_i \times \\
		m_i\:\:{c_i\times}^T  & m_i\mathrm{1} \end{bmatrix},
\end{eqnarray}
which is primarily obtained from the momentum equation of \eqref{eq:momentum}. $\bar{I}_{C_i}$ is the inertia matrix taken w.r.t. the center of mass and $m_i$ is the mass for link $i$. $c_i$ is the center of mass location of the same link w.r.t. the joint axis $O_i$. Accordingly, the kinetic energy of the $i^{th}$ link is given by using \eqref{eq:bodyvelocity}. Substituting for $v_{O_i}$ in \eqref{eq:bodyvelocity} results in:
\begin{eqnarray}
\label{eq:energysubs}
T_i= \frac{1}{2} \dq^T J_{i}^T \left [ \begin{array}{cc} \nu_i & 0 \\ 0 & \nu_i \end{array} \right ]^T I_{O_i} \left [ \begin{array}{cc} \nu_i & 0 \\ 0 & \nu_i \end{array} \right ] J_{i} \dq,
\end{eqnarray}
where $\nu_i$ denotes the rotation of the $i^{th}$ link w.r.t. the joint axis. The inertia matrix, $D(\q)$, can thus be expressed from the total energy of the $b$ bodies:
\begin{eqnarray}
 \label{eq:totalkineticenergy}
 T &=& \sum\limits_{i=1}^b T_i  \\
 &=& \frac{1}{2} \dq^T \underbrace{\left ( \sum\limits_{i=1}^b J_{i}^T \left [ \begin{array}{cc} \nu_i & 0 \\ 0 & \nu_i \end{array} \right ]^T I_{O_i} \left [ \begin{array}{cc} \nu_i & 0 \\ 0 & \nu_i \end{array} \right ] J_{i} \right )}_{D(\q)} \dq \nonumber,
\end{eqnarray}
where the inertia tensor $I_{O_i}$ is obtained from \eqref{eq:inertiatensorspecific}.
Note that the Jacobian $J_i$ is purely a function of the joint angles and the link lengths. Therefore, assuming that the distances between the joint axes are known (which are easy to measure), all the other terms, namely, center of mass position, inertia, masses are in the spatial inertia tensor $I_{O_i}$. This fact will be utilized in later sections to compute the regressor.
%


$C(\q,\dq)$ can also be derived as a linear function of the same elements of the tensor, by utilizing the Christoffel symbols: 
$\Gamma_{ijk}$ is obtained from the inertia matrix $D(\q)$:
\begin{eqnarray}
\label{eq:coriolisdetails}
\Gamma_{ijk} &=& \frac{1}{2} \left ( \frac{\partial{D_{ij}(\q)}}{\q_k}+\frac{\partial{D_{ik}(\q)}}{\q_j}-\frac{\partial{D_{kj}(\q)}}{\q_i} \right )
\nonumber \\
 C_i (\q,\dq)  \dq &=& \sum\limits_{j,k=1}^b {\Gamma}_{ijk} {\dq}_{j} {\dq}_{k}.
\end{eqnarray}

The potential energy function $V(\q)$ can be computed as sum of the potential energies of the individual links:
\begin{eqnarray}
 \label{eq:potentialenergy}
 V(\q) = \sum\limits_{i=1}^b V_i(\q) =\sum\limits_{i=1}^b m_i g h_i(\q),
\end{eqnarray}
where $m_i$ is the mass of the individual links, $g$ is the gravity and $h_i$ is the vertical position of the center of mass for each link. Therefore, $h_i$ is the sum of heights of the $i^{th}$ joint axis, $h_{O_i}$, and the vertical height of the CoM w.r.t. the joint axis:
\begin{eqnarray}
 \label{eq:heights}
 h_i(\q) = h_{O_i}(\q) + h_{C_i} (\q),
\end{eqnarray}
since the link lengths are assumed to be known, $h_{O_i}$ in \eqref{eq:heights} is known. The height of the CoM $h_{C_i}$ is the dot product of the center of mass location $c_i=[c_{i,1},c_{i,2},c_{i,3}]^T$ rotated by a transformation matrix, $\nu$, and the vertical axis. Assuming that $z$ axis is along the vertical axis, we have:
\begin{eqnarray}
 \label{eq:heightcom}
 h_{C_i}(\q) = \left[\begin{array}{c} 0 \\ 0 \\ 1 \end{array} \right ]^T \nu(\q)  c_i.
\end{eqnarray}

\eqref{eq:heightcom} and \eqref{eq:heights} can be substituted in \eqref{eq:potentialenergy} to obtain:
\begin{eqnarray}
 \label{eq:potentialmod}
 V(\q) = \sum\limits_{i=1}^n m_i g h_{O_i}(\q) + \sum\limits_{i=1}^b g  \underbrace{\left[\begin{array}{c} 0 \\ 0 \\ 1 \end{array} \right ]^T \nu(\q) }_{\kappa}  m_i c_i,
\end{eqnarray}
where the unknown parameters are linear in the expression and are already present in the inertia tensor $I_{O}$. Therefore, in perspective, all the unknown parameters are effectively collected in $I_O$, which motivates the path that this paper takes to compute the regressor for any general $n$-DOF $b$-body system.

\section{The Regressor and the Parameters }
\label{sec:reg}

Since the parameters are not perfectly known, the equation of motion, \eqref{eq:eom} computed with the given set of parameters will be henceforth have {\Large $\hat{}$} over the symbols. Therefore, $D_a, C_a, G_a$ are the actual inertia, motor inertia, Coriolis and gravity matrices of the robot, and ${\hat D}, {\hat C}, {\hat G}$ are the assumed inertia, Coriolis and gravity matrices of the robot. 

Consider the equation of motion of an n-link robot which is obtained from the Lagrangian \eqref{eq:Lag} and is restated here as:
\begin{eqnarray}
 \label{eq:eomidentification}
K(\q,\dq,\ddq)=B \mathrm{U} ,
\end{eqnarray}
where $K={ D}(\q) \ddot{\q} + { C}(\q,\dot{\q}) \dot{\q} + { G}(\q)$ obtained from \eqref{eq:eom}. 

It is a well known fact that the parameters of a robot, like the inertia, masses, position of center of mass are affine in \eqref{eq:eomidentification} (see \cite{spong2006robot}). Therefore, it is possible to write \eqref{eq:eom} in the form:
\begin{eqnarray}
 \label{eq:separatedeom}
 K=\mathrm{Y}(\q,\dq,\ddq) \Theta = \mathrm{U},
\end{eqnarray}
where $\mathrm{Y}(\q,\dq,\ddq)$ is called the regressor in \cite{spong2006robot}, and $\Theta \in \R^{n_P}$ is called the set of base inertial elements (parameters). It is important to note that $\Theta$ need not be unique, and if the set is unique, then it is called the Base Parameter Set \cite{ITRA:bpsmayeda}. We can write, $\Theta=[\theta_1,\theta_2,\theta_3,\dots]^T$, where each $\theta_i$ is a function of the unknown parameters of the robot. $n_P$ is the size of the parameter set. Accordingly, $K_a=\mathrm{Y}(\q,\dq,\ddq) \Theta_a$, and $\hat{K}=\mathrm{Y}(\q,\dq,\ddq) \hat{\Theta}$, where $\Theta_a$ is the actual set of base inertial parameters, and $\hat{\Theta}$ is the assumed set of base inertial parameters. 

\newsec{Determining the base inertial parameters ($\Theta$).} Since we know the inertia of link $i$, $\bar{I}_C \in \R^{3 \times 3}$ is symmetric, and the square of a skew-symmetric matrix results in a symmetric matrix, $\bar{I}_{C_i}+ {c_i\times} {c_i\times}^T$ is symmetric. We can assign the inertial parameters of the $i^{th}$ link $\Theta_i=[\theta_{i,1},\theta_{i,2},\dots]$ to the elements of the spatial inertia tensor $I_{O_i}$ given in \eqref{eq:inertiatensorspecific} in the following manner:
\begin{eqnarray}
\label{eq:assigningtensor}
I_{O_i} = \left [ \begin{array}{cccccc} \theta_{i,1} & \theta_{i,2} & \theta_{i,3} & 0 & -\theta_{i,4} & \theta_{i,5} \\
                                        \theta_{i,2} & \theta_{i,6} & \theta_{i,7} & \theta_{i,4} & 0 & -\theta_{i,8} \\
										\theta_{i,3} & \theta_{i,7} & \theta_{i,9} & -\theta_{i,5} & \theta_{i,8} & 0 \\
										0 & \theta_{i,4} & -\theta_{i,5} & \theta_{i,10} & 0 & 0 \\
										-\theta_{i,4} & 0 & \theta_{i,8} & 0 & \theta_{i,10} & 0 \\
										\theta_{i,5} & -\theta_{i,8} & 0 & 0 & 0 & \theta_{i,10} \end{array} \right ] \nonumber ,
\end{eqnarray}
where $10$ parameters are obtained for each link. This is similar to how the parameters were categorized in \cite{khosla1989categorization}, which uses the newton-euler method directly. The major difference is that \eqref{eq:assigningtensor} is obtained from the tensors and are used directly in online identification and control, which gets tedious without SVA. Obtaining $D(\q)$ and $C(\q,\dq)$ from $I_{O_i}$ are straightforward from \eqref{eq:totalkineticenergy} and \eqref{eq:coriolisdetails}. Consider the gravity vector, which is the partial derivative w.r.t $\q$ of the potential energy, $V(\q)$:
\begin{eqnarray}
G(\q)   = \sum\limits_{i=1}^n g \left ( \frac{\partial{h_{O_i}}}{\partial{\q}} {m_i} +  \frac{\partial{\kappa}}{\partial{\q}} {c_i m_i} \right )
\end{eqnarray}
Therefore, even the gravity vector $G(\q)$ is a linear function of the inertial parameters present in the tensor.

Comparing with \cite{atkeson1986estimation} and \cite{dhparam}, the inertial parameters chosen were different than the one chosen here. Specifically, the inertial parameters in $\bar{I}_C$ form the base inertial parameters; whereas here the parameters in the matrix $\bar{I}_{C_i}+ {c_i\times} {c_i\times}^T$ make the unknown parameter set. Besides, it is also possible to directly compute the regressor while evaluating the dynamics of the robot, by just picking the coefficient of every parameter $\theta_i$ one by one. In fact, this becomes the basis for a very simple algorithm shown in \algref{alg:regressor}. For a robot having $b$ rigid bodies, the number of unknown parameters will be $10b$.


It is shown in \cite{featherstone:RBDL} that it is possible to compute inverse dynamics of n-DOF robot in $O(n)$. This is achieved by deploying Newton-Euler recursive method. This is a standard technique used in numerical computation of the dynamics, which was used initially in the 1980's. Therefore, assuming that it is possible to compute the eom
in just $n$ recursive iterations, we propose \algref{alg:regressor} which calls in the current state and acceleration of the robot and computes the regressor from the data.
 \begin{algorithm}[t]
 \caption{Regressor Pseudocode}
 \label{alg:regressor}
 \begin{algorithmic}
 \FOR {$i=1$ \TO $n$}
 \FOR {$j=1$ \TO $10$}
 \STATE {$\theta_{i,j}=0$}
 \ENDFOR
 \STATE {Update $I_{O_i}$ with the value $\theta_{i,j}$}
 \ENDFOR
%
 \FOR {$j=1$  \TO $n$} 
 \FOR {$j=1$ \TO $10$}
 \STATE {$\theta_{i,j}=1$}
 \STATE {Update $I_{O_i}$ with the value $\theta_{i,j}$}
 \STATE {$\mathrm{Y}_{i+j-1} = D(\q)\ddq+C(\q,\dq)\dq+G(\q)$}
 \STATE {$\theta_{i,j}=0$}
 \STATE {Update $I_{O_i}$ with the value $\theta_{i,j}$}
 \ENDFOR
 \ENDFOR
 \end{algorithmic}
 \end{algorithm}
The number of iterations for evaluating the regressor is $10b$, which is proportional to the number of rigid bodies present in the robot. And the maximum number of degrees of freedom for each rigid body is $6$, which implies that $b \leq n \leq 6b$. Accordingly in each iteration the equation of motion is computed which takes $n$ iterations, and the resulting algorithmic complexity will be $O(n^2)$.

\section{Estimation and Control}
\label{sec:estctrl}

The regressor of the previous section has two main applications 1) to facilitate the identification of unknown model parameters and 2) to enable straightforward calculation of computed-torque controllers.  This section states these problems in the general case, with specific examples to follow in Section \ref{sec:casestudy}. 

\newsec{Parameter Identification.}
The problem of parameter identification can be stated as follows: Suppose we are given a fully actuated robotic linkage with $b$ rigid bodies, $n$ degrees of freedom and unknown inertial parameters $\Theta_a \in R^{10\Ndof}$.  Given $s$ vectors of torque $\mathrm{U}=[u_1, u_2, \dots, u_{n}]^T$, generalized configuration $q=[q_1, q_2,$ $ \dots, q_{n} ]^T$, generalized velocity data $\dq = [ \dot{q}_1, \dot{q}_2, \dots, \dot{q}_{n} ]^T$, and generalized acceleration $\ddq=[ \ddot{q}_1, \ddot{q}_2, \dots, \ddot{q}_{n} ]^T$, choose model parameters $\hat{\Theta} \in R^{10\Ndof}$ such that\footnote{If necessary, the statement of the parameter identification problem \eqref{eqn:parameter_ID} can be modified to include a requirement that $\Theta > 0$.  The resulting problem will, in general, no longer have a closed form solution and could instead be solved using constrained quadratic programming.}

\begin{eqnarray} \label{eqn:parameter_ID}
\hat{\Theta} = \underset{{\Theta} \in \R^{10\Ndof}}{\operatorname{argmin}} & \:\:\:\:\: \| \U - \Y \Theta \|_2 
\end{eqnarray}
where $\U$ is the collection of torque vector inputs and $\Y$ is the collection of regressor matrices for $s$ samples of angle, velocity and acceleration data. $\U$ and $\Y$ are given as:
\begin{equation*}
\begin{array}{cc}
\U = \left[ \begin{array}{c} \mathrm{U}[1] \\ \mathrm{U}[2] \\ \vdots \\ \mathrm{U}[s] \end{array} \right], 
\Y = \left[ \begin{array}{c} \mathrm{Y}(q[1],\dot{q}[1],\ddot{q}[1]) \\ \mathrm{Y}(q[2],\dot{q}[2],\ddot{q}[2]) \\ \vdots \\ \mathrm{Y}(q[s],\dot{q}[s],\ddot{q}[s]) \end{array} \right]
\end{array}
\end{equation*}
and
\begin{equation*}
\begin{array}{ccc}
q[i] = \left[ \begin{array}{c} q_1[i] \\ q_2[i] \\ \vdots \\ q_n[i]\end{array}\right] & \dot{q}[i] = \left[ \begin{array}{c} \dot{q}_1[i] \\ \dot{q}_2[i] \\ \vdots \\ \dot{q}_n[i] \end{array}\right] & \ddot{q}[i] = \left[ \begin{array}{c} \ddot{q}_1[i] \\ \ddot{q}_2[i] \\ \vdots \\ \ddot{q}_n[i] \end{array}\right]. 
\end{array}
\end{equation*}
A vector $\hat{\Theta}$ that minimizes $\| \U - \Y \hat{\Theta} \|_2$ can be found using the Moore-Penrose pseudoinverse:
\begin{equation} \label{eqn:theta_hat}
\hat{\Theta} = \text{pinv}(\Y)\U.
\end{equation}
An estimated (or modeled) set of computed torques $\Uhat$ can be calculated using the parameter vector $\hat{\Theta}$,
\begin{equation}
\label{eq:torqueestimate}
\Uhat = \mathrm{Y}_C\hat{\Theta}.
\end{equation}
Referring to the definition of $\Uhat$ above, a set of estimated torques can be found for each actuator.  These estimates will be denoted $[\hat{u}_1, \hat{u}_2, \dots, \hat{u}_n]$.  The coefficient of determination, $R^2$, can be used to describe the similarity between $\U$ and $\Uhat$ (or equivalently between $[u_1, u_2, \dots, u_n]$ and $[\hat{u}_1, \hat{u}_2, \dots, \hat{u}_n]$),
\begin{equation}
\begin{array}{l}
e = \U - \Y\hat{\Theta} \\
R^2 = 1-(e^Te)/(\U^T\U).
\end{array}
\end{equation}

\newsec{Parameter Identification with an Initial Guess.}
An initial guess or nominal parameter value can be incorporated into the parameter identification problem by modifying the cost function in \eqref{eqn:parameter_ID}.

\begin{align} \label{eqn:parameter_ID_extended}
\hat{\Theta} = \underset{{\Theta} \in \R^{10\Ndof}}{\operatorname{argmin}} & \:\:\:\:\: \alpha \| \U - \Y \Theta \|_2 + (1-\alpha) \|\Theta - \Theta_0\|_2
\end{align}
where $\Theta_0 \in R^{10\Ndof}$ is a nominal parameter vector, and $\alpha \in (0,1)$ is a factor that can be used to vary the relative effects of the problem data ($\U$ and $\Y$) and the initial guess ($\Theta_0$).  A least squares solution to \eqref{eqn:parameter_ID_extended} can be found that is similar in structure to the solution of \eqref{eqn:parameter_ID},

\begin{equation} \label{eqn:theta_hat_extended}
\hat{\Theta} = \text{pinv}(\widetilde{\Y})\widetilde{\U}
\end{equation}
where
\begin{equation}
\label{eq:pinv_matrices}
\widetilde{\U} = \left[ \begin{array}{c} \alpha \, \U \\ (1-\alpha) \,\Theta_0 \end{array}\right], \,\,  
\widetilde{\Y} = \left[ \begin{array}{c} \alpha \, \Y \\ (1-\alpha) \,{\bf I}_{10 \Ndof} \end{array} \right].
\end{equation}

 \newsec{Rank Properties of the Regressor.}
 It is not necessary to consider all of the samples to compute the parameters because samples leading to low eigen values of the matrix $\Y$ form a computation burden.
 Furthermore, it is possible to obtain $p^*$ samples which give the parameter estimate $\hat{\Theta}=\Theta^*$ from the optimization problem \eqref{eqn:parameter_ID} such that:
 \begin{eqnarray}
\label{eq:equalitypistar}
\Y(\q,\dq,\ddq) {\Theta}^*=\Y(\q,\dq,\ddq) {\Theta}_a
\end{eqnarray}

Once a parameter vector is identified, the regressor can be used in a computed torque controller to bring about a desired acceleration of the joints.  Let $\ddot{q}_{cmd} \in R^{n}$ represent a desired vector of joint accelerations, then a unique controller to induce these accelerations is given by:
\begin{equation}
\label{eq:computedtorque}
\mathrm{U}_{cmd} = \mathrm{Y}(q,\dot{q},\ddot{q}_{cmd})\hat{\Theta}.
\end{equation}
The following lemma will introduce the relationship between the desired and actual acceleration of the robot.
\begin{lemma} {\it
 For the fully actuated robot, i.e., $B= I_{n\times n}$, if $\hat{\Theta}=\Theta^*$ is evaluated from the optimization problem \eqref{eqn:parameter_ID} with $p^*$ samples, and if the control law used is \eqref{eq:computedtorque}, then $\ddq=\ddot{q}_{cmd}$.}
\end{lemma}

Note that the inertial parameters obtained from above will not yield true parameters of the robot, but will give the same value for the computed torque as described by Lemma 1. This property will be used in implementing online model based controllers and eliminate the identification of true parameters of the robot.

\section{Experimental Results}
\label{sec:casestudy}
This section presents three experimental studies which use the regressor of \secref{sec:reg} to solve problems of identification and the control of identified systems. The first application will be in offline identification of a planar double pendulum, the second in offline identification of a 3D robotic leg of the Valkyrie robot (see \figref{fig:valk}), and the third in online identification and control of a 3D robotic arm of the Valkyrie robot (see \figref{fig:valk}).

\begin{figure}[ht!] 
  \centering
  \begin{subfigure}[b]{.3\columnwidth}
  \includegraphics[width=\columnwidth]{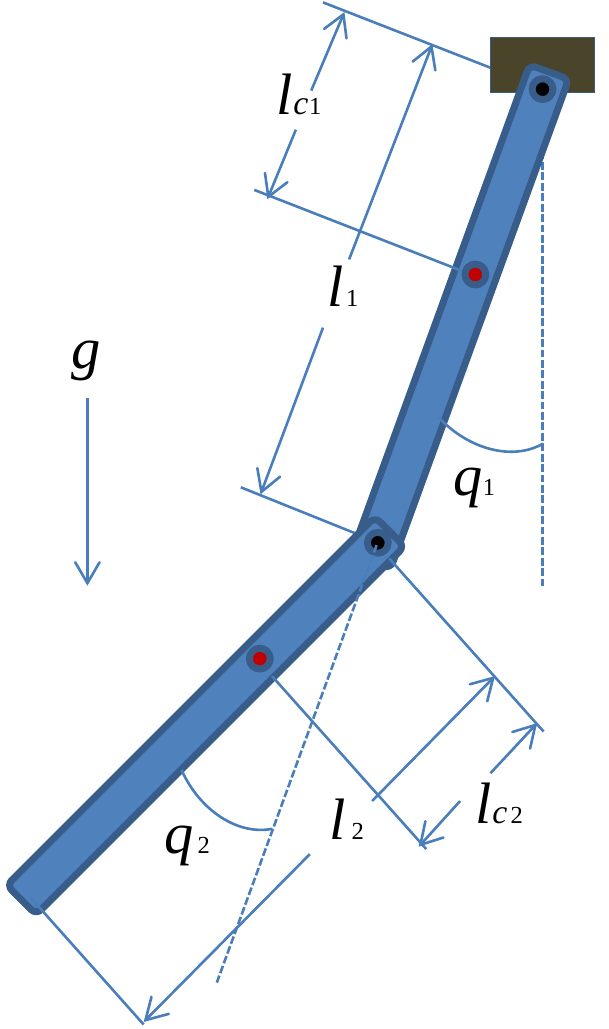}
  \caption{Illustration of the double pendulum.} 
  \label{fig:doublependulum_img}
  \end{subfigure}\hspace{1em}
  \begin{subfigure}[b]{.52\columnwidth}
  \includegraphics[width=\columnwidth]{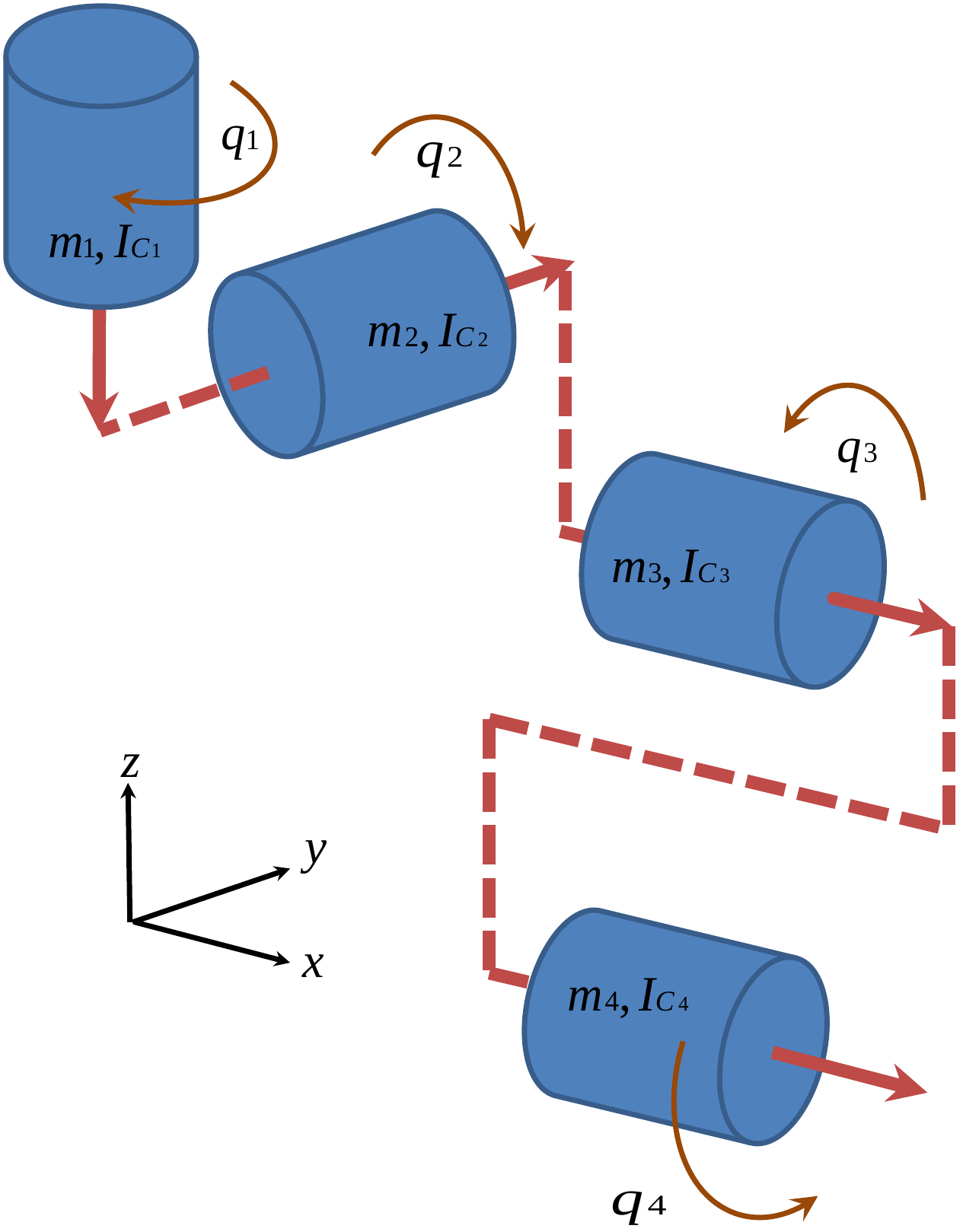}
  \caption{Figure showing the 4-link leg.} 
  \label{fig:robotleg}
  \end{subfigure} \\[1em]
  
  \begin{subfigure}[b]{.75\columnwidth}
  \includegraphics[width=\columnwidth]{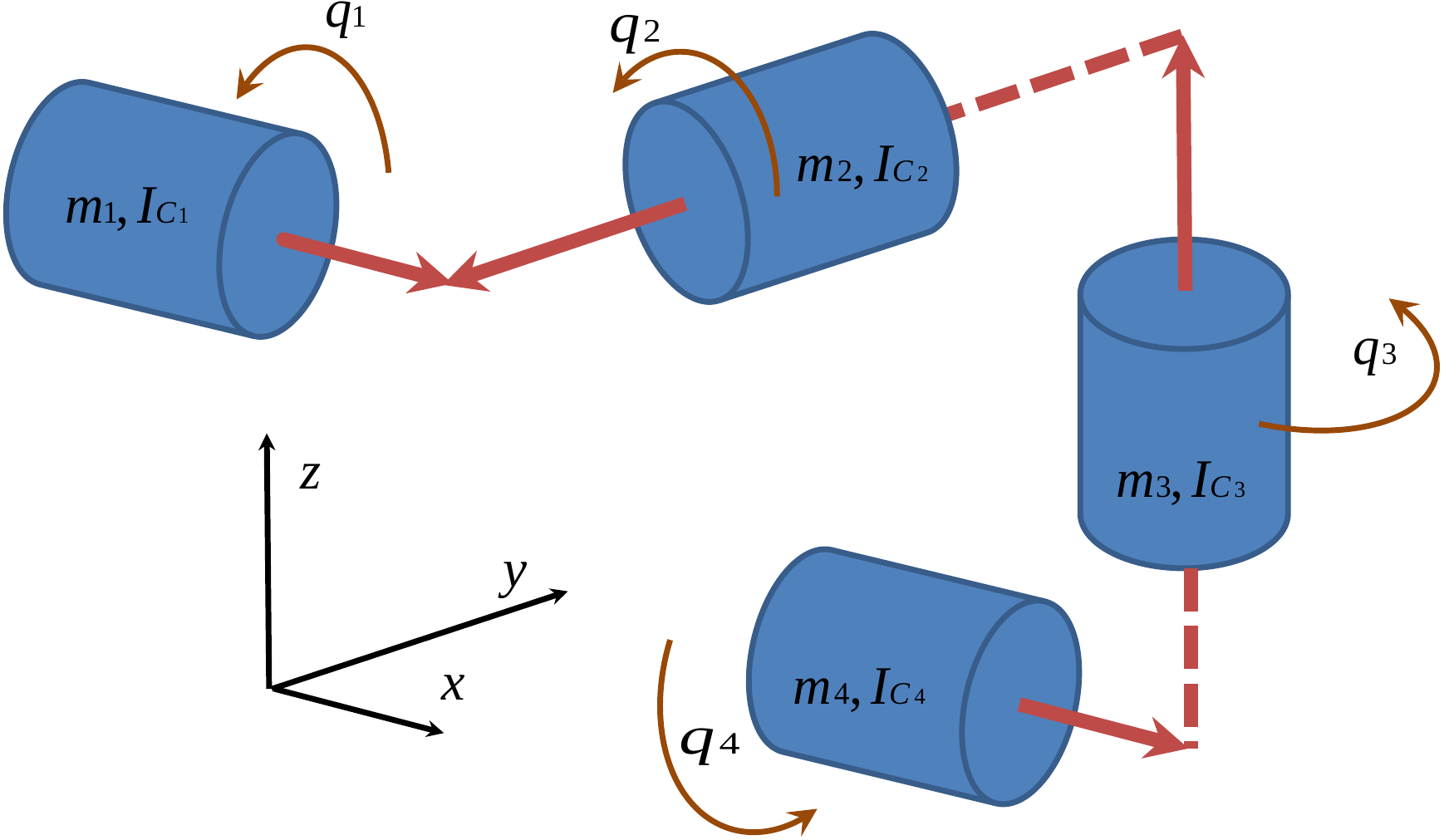}
  \caption{Figure showing the robot arm.} 
  \label{fig:robotarm}
  \end{subfigure}\\

  \caption{Robot models considered in this paper.}
  \label{fig:robots}
\end{figure}

\begin{figure}[t!]
  \centering
  \begin{subfigure}[b]{\columnwidth}
    \includegraphics[width=0.99\columnwidth]{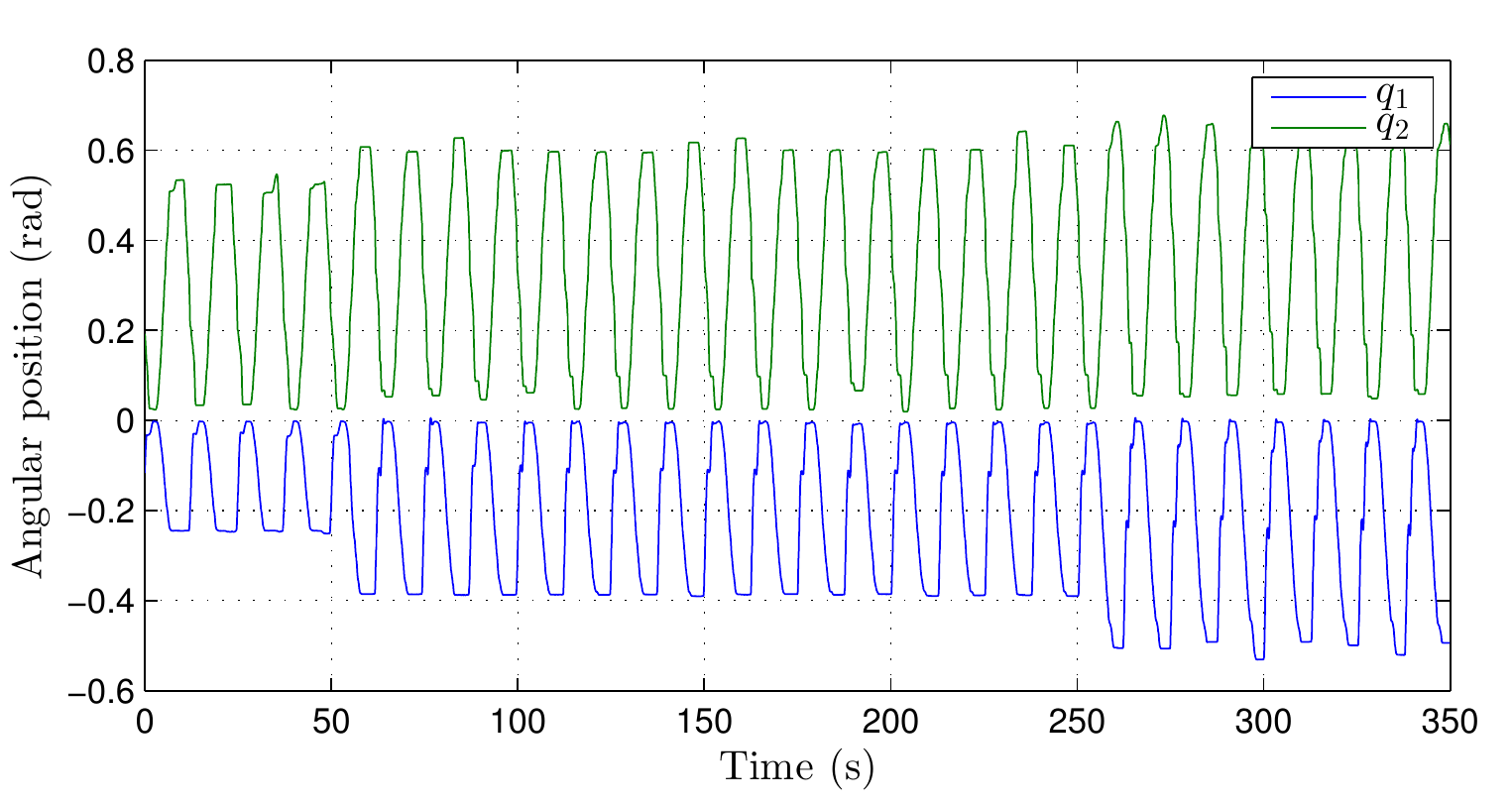}
    \vspace{-2em}
  \end{subfigure}

  \begin{subfigure}[b]{\columnwidth}
    \includegraphics[width=0.99\columnwidth]{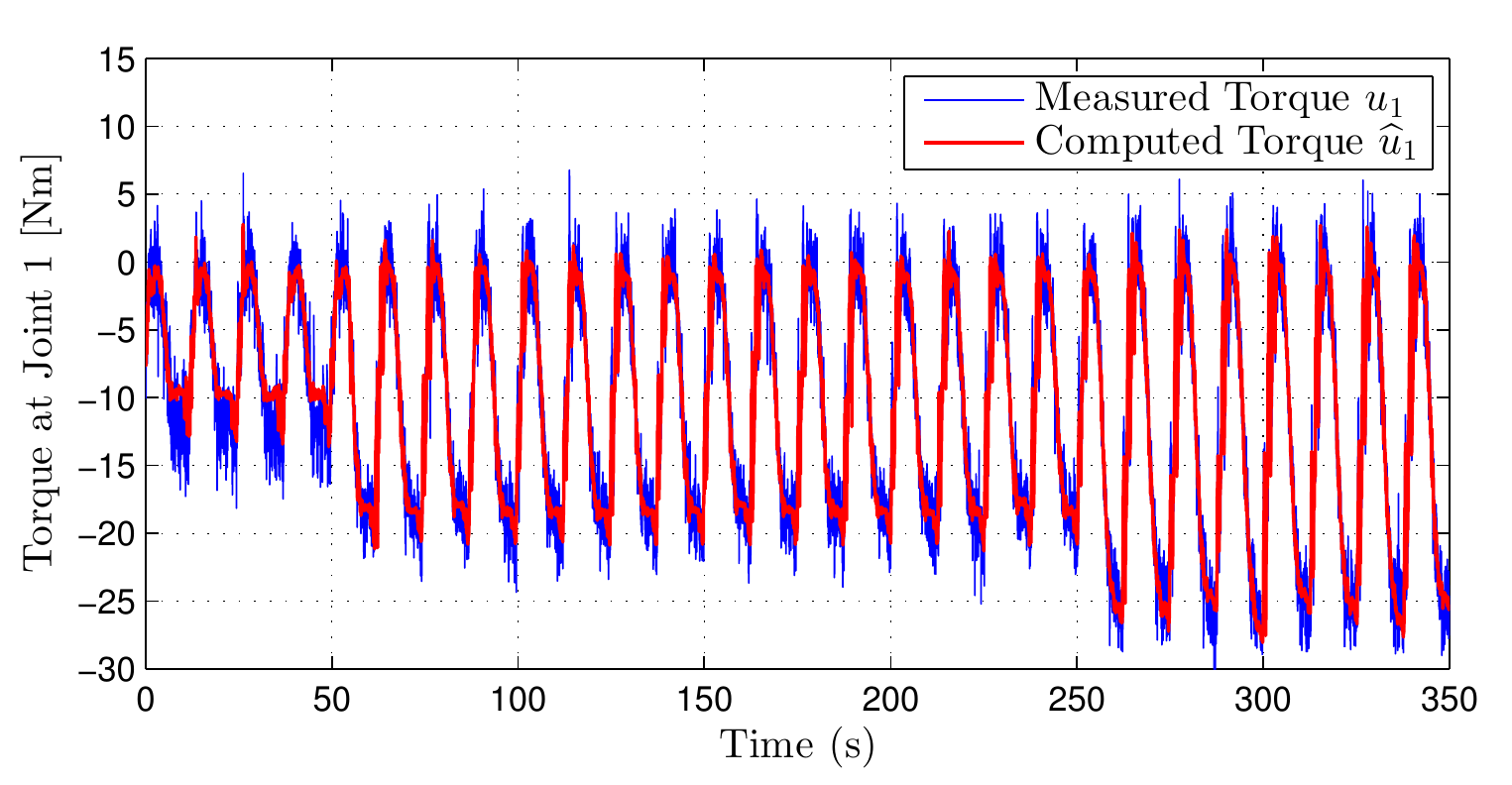}
    \vspace{-2em}
  \end{subfigure}

  \begin{subfigure}[b]{\columnwidth}
    \includegraphics[width=0.99\columnwidth]{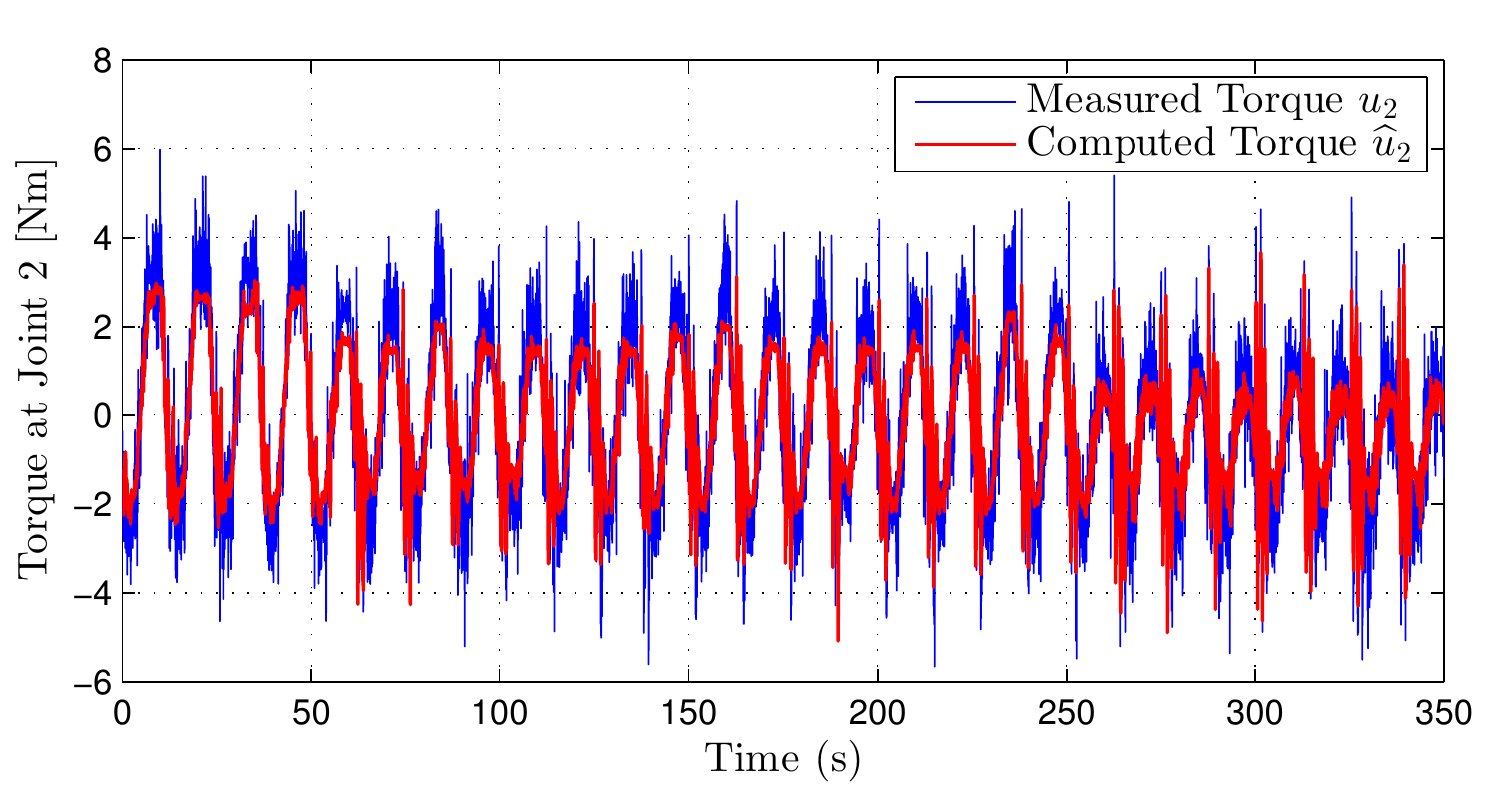}
    \vspace{-1em}
  \end{subfigure}
    \caption{\small {\bf Top:} Joint trajectories as a function of time experienced by the experimental setup of Figure \ref{fig:doublependulum_img}.  {\bf Middle, Bottom:} A comparison between experimentally measured torque vectors $u_1$, $u_2$ and the computed torque values $\hat{u}_1, \hat{u}_2$ corresponding to the identified system.}
  \label{fig:doublependulum_fit}
\end{figure}

\begin{figure}[t!]
  \centering
  \includegraphics[width=0.95\columnwidth]{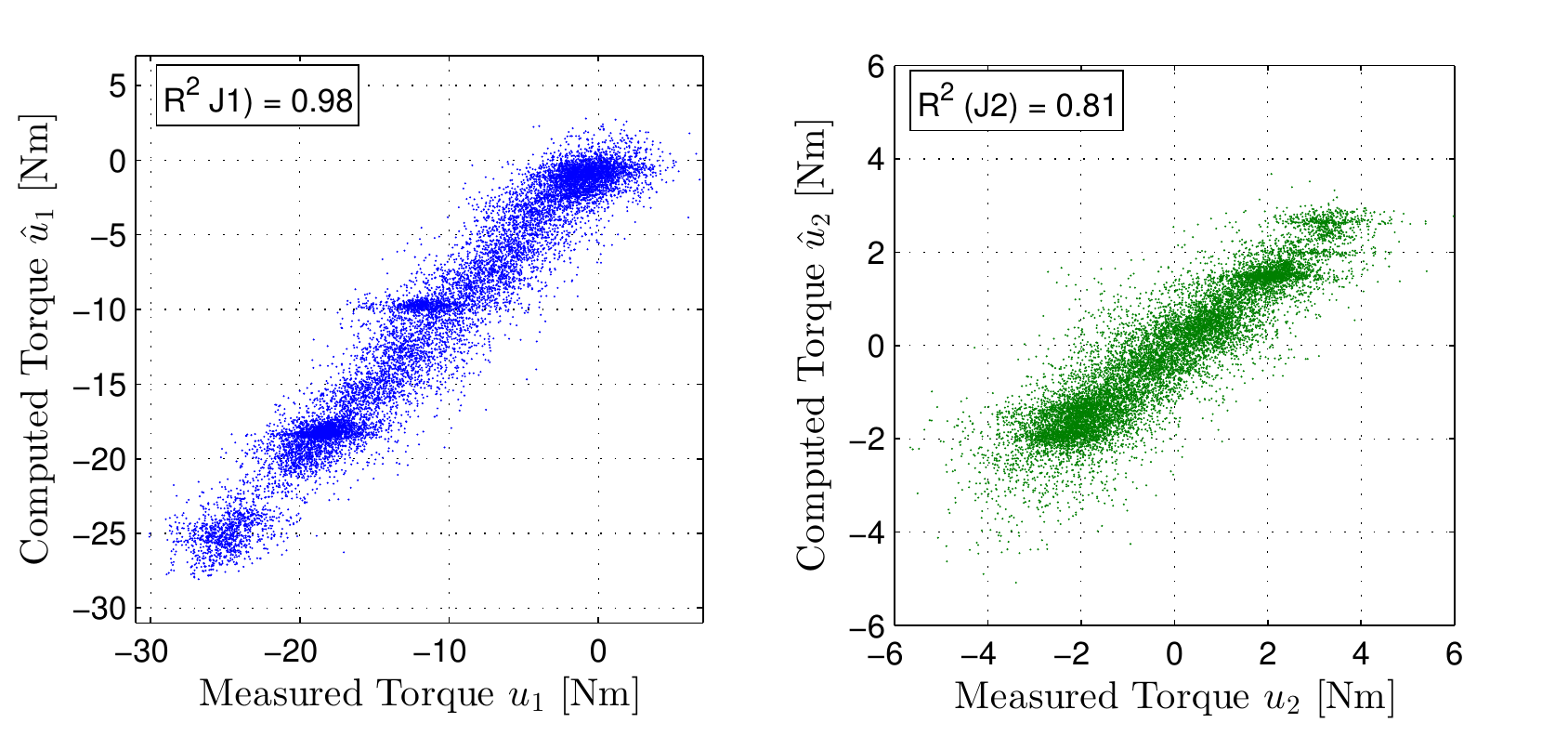}
  \vspace{-1em}
  \caption{\small Plotting measured torque vs. computed torque provides an illustration of quality of fit.  The coefficient of determination, $R^2$, is shown for each joint.}
  \label{fig:doublependulum_rsquared}
\end{figure}

\newsec{Offline Identification of a Planar Double Pendulum.}
Consider a double pendulum, such as the one illustrated in \figref{fig:doublependulum_img} with inertia tensors for both link 1 and 2 which are computed from \eqref{eq:inertiatensorspecific}.
Note that since the pendulum considered is planar, the number of parameters used from the inertia tensor effectively being used is $3$. 

\figref{fig:doublependulum_fit} shows the outcome of an experiment where angles, velocities, accelerations, and torques were measured at each controller instant while the linkage executed a sinusoidal motion at each joint.  An estimate of the parameter vector $\Theta$ was made without reference to an initial guess, and thus the estimation scheme of \eqref{eqn:theta_hat} was used. The tracking plots of \figref{fig:doublependulum_fit} and the scatter plots of \figref{fig:doublependulum_rsquared} illustrate that the offline identification procedure (using all collected data for a single parameter fit) produced a parameter set leading to very high correlation.  As noted in \secref{sec:reg}, the regressor $\Y$ will not necessarily be full-rank, and thus any estimate $\hat{\Theta}$ should not be expected to converge to the actual parameters $\Theta_a$, even for very large data sets.

\begin{figure}[t!] \centering
  \includegraphics[width=\columnwidth]{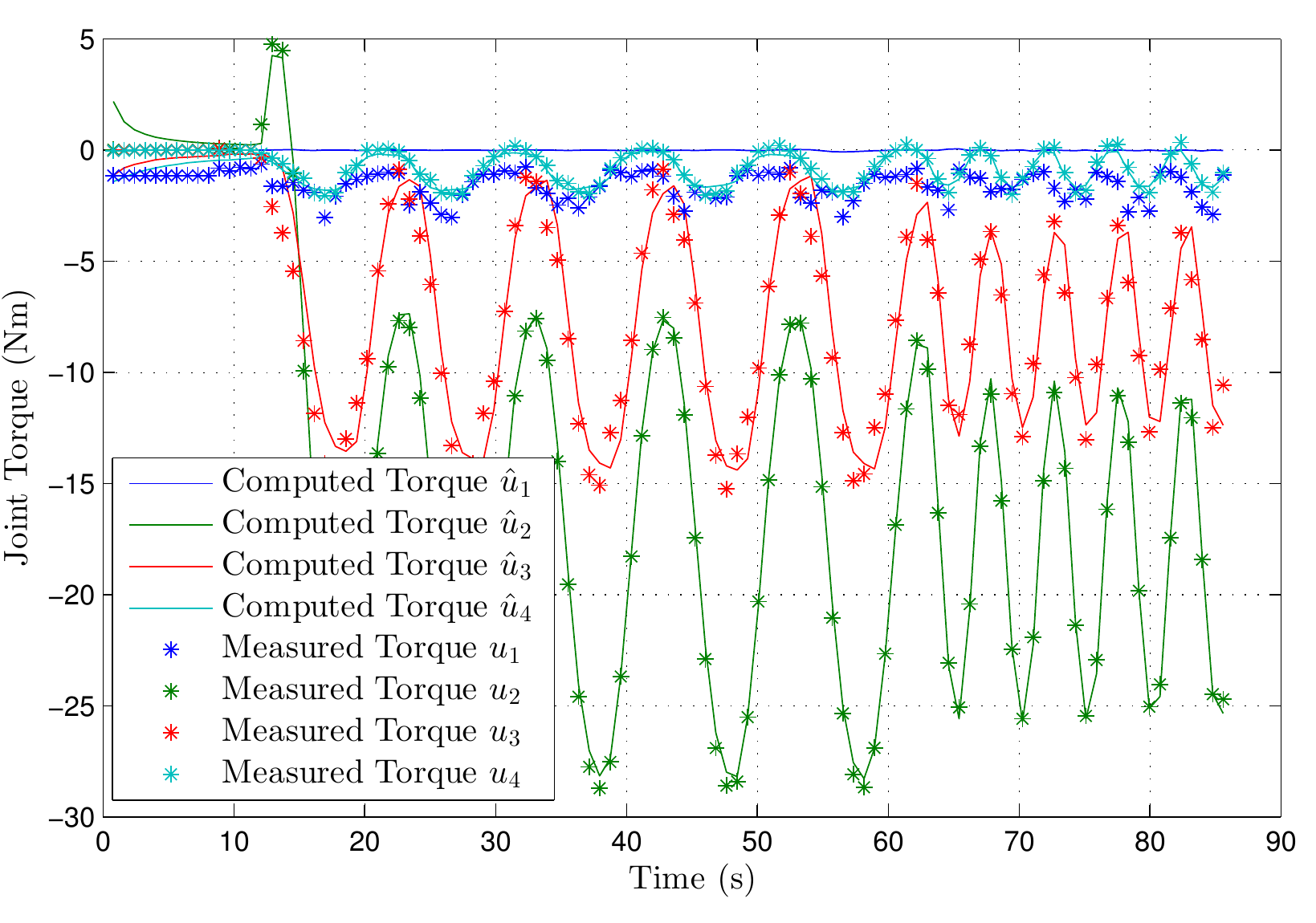}
  \includegraphics[width=\columnwidth]{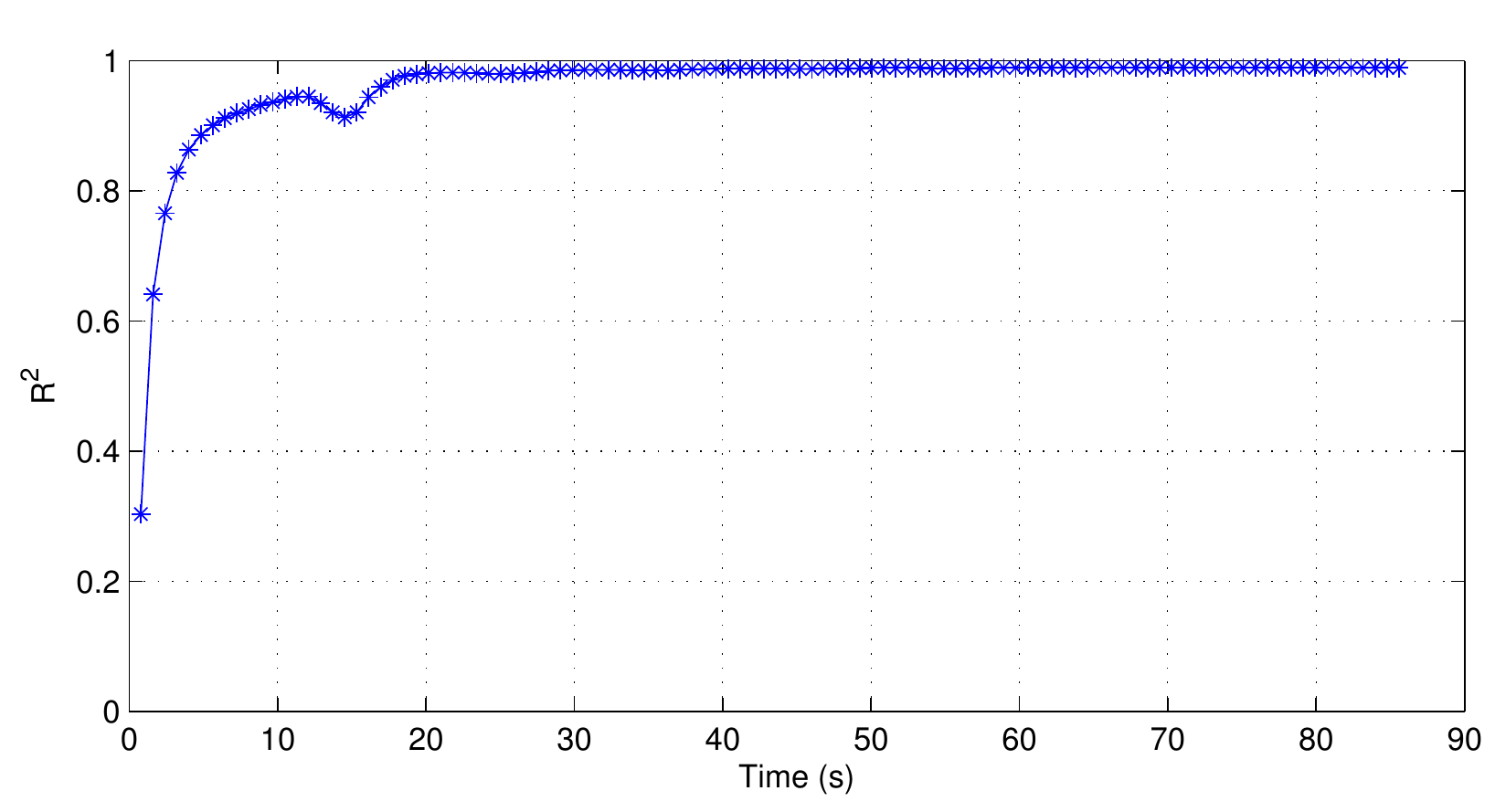}
  \caption{\small The plots above show the quality of fit improves during an online identification exercise of the robot leg pictured in \figref{fig:robotleg}. The top figure shows the instantaneous computed torque as a function of time, along with the measured data that were used for the parameter fit. Note that both the $R^2$ plot at bottom and the torque tracking plot at top both indicate that the output of the identification algorithm leads to noticeably higher quality of fit as the length of the available data vectors increases.} 
  \label{fig:leg_rsquared}
\end{figure}

\newsec{Offline Identification of a 4-Link Robotic Leg.}  4-DOF leg shown in \figref{fig:robots} and \figref{fig:valk} was identified offline. In the experiment a series of position, velocity, acceleration, and torques were measured at each of the four joints of the robot.  An online estimate was made of the data, where the parameter identification was updated at a lower, decimated rate.  The fit algorithm used here for online identification is very similar to \eqref{eqn:theta_hat_extended} used in offline identification of the double pendulum.  The two differences for online application are that 1) the data vectors $\U$ and $\Y$ grow throughout the experiment, and as such the quality of fit improves the longer the experiment is run, and 2) the initial parameter guess $\Theta_0$ is used to ensure bounded behavior at the beginning of the experiment.  Results of the fit are shown in \figref{fig:leg_rsquared}.  

\begin{figure}[t!]
  \centering
  
  \begin{subfigure}[b]{1.0\columnwidth}
    \includegraphics[width=\columnwidth]{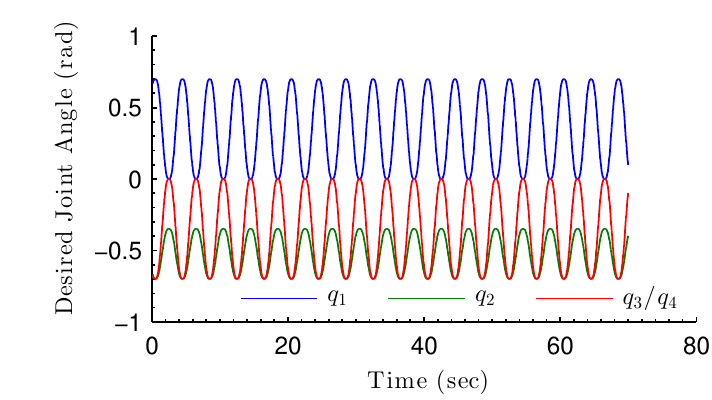}
    \caption{%
The arm was made to exhibit a tick tock behavior where the joints followed the trajectories shown.}
    \label{fig:3d-arm-desired}
  \end{subfigure}
  
  \begin{subfigure}[b]{1.0\columnwidth}
    \includegraphics[width=\columnwidth]{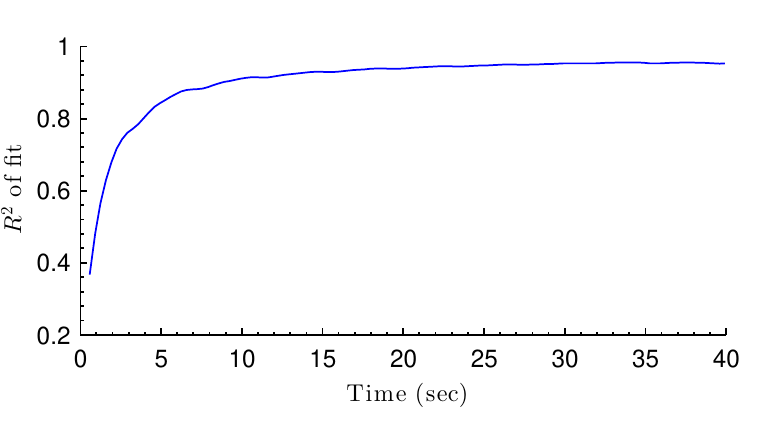}\vspace{-1em}
    \caption{%
The coefficient of determination, $R^{2}$, can be seen to converge to right below unity.
With perfect sensing and no dependence on an incorrect initial guess, this quantity will converge precisely to unity.
Due to noise, the evolution is non-monotonic.}
    \label{fig:3d-arm-r2}
  \end{subfigure}

  \begin{subfigure}[b]{1.0\columnwidth}
    \includegraphics[width=\columnwidth]{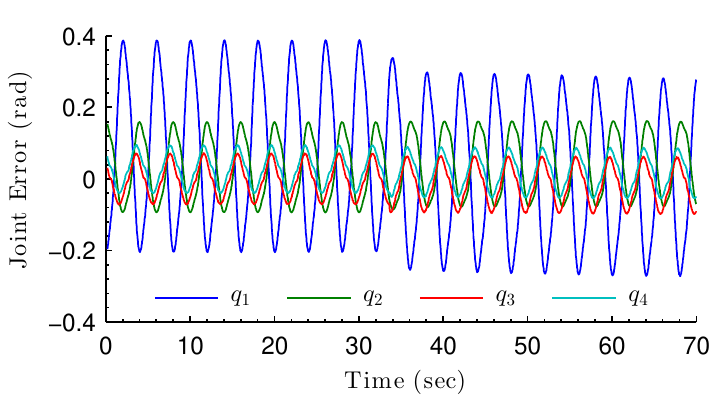}
    \vspace{-1em}
    \caption{%
The peak-to-peak amplitudes of the joint errors appear essentially unchanged when the controller switches to using the identified model around 35 seconds, but the bias is clearly shifted closer to zero.}
    \label{fig:3d-arm-errors}
  \end{subfigure}
  \caption{Experimental results for a 4-DOF robotic arm.}
\end{figure}

\begin{figure*}
\centering
\includegraphics[width=0.7\textwidth]{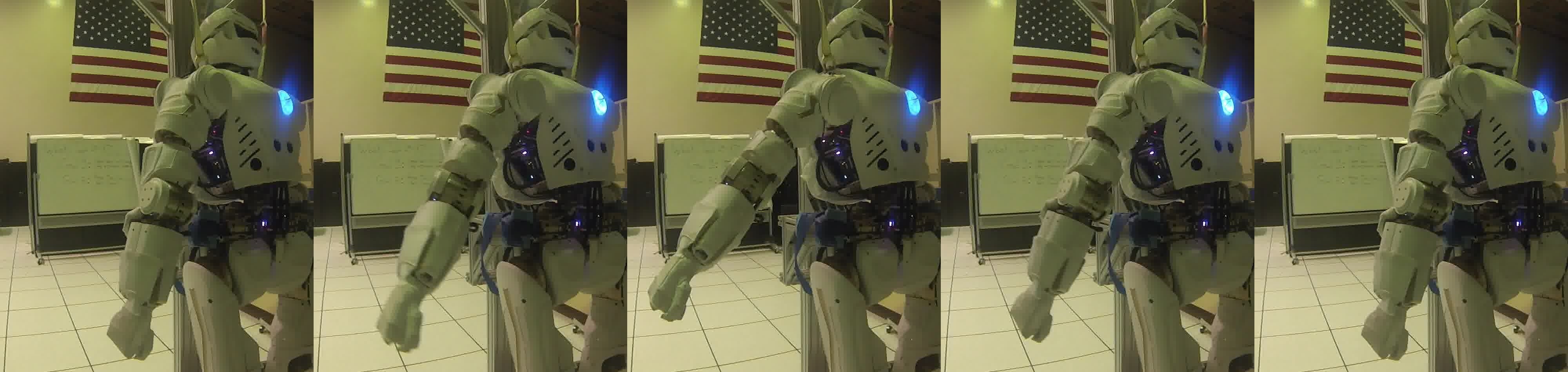} 
\caption{Figure showing the tile of 4-DOF arm swing experiment used for identification.}
\label{fig:robotarmvalk}
\end{figure*}

\newsec{Online Identification of a 4-Link Robotic Arm.}
%
%
Online identification was conducted on a robotic arm with four degrees of freedom (see \figref{fig:robotarm} and \figref{fig:valk}).
To draw an analogue to human physiology, the joints $q_1$, $q_2$, $q_3$, and $q_4$ can be thought of as the shoulder extensor, the shoulder adductor/abductor, the upper arm pronator /supinator, and the elbow extensor, respectively.
The joint angles and even velocities, which are necessary for computing the regressor, can be measured from encoders but accelerations must also be obtained and a common procedure is to filter accelerations -- in this experiment we used an exponential moving average filter.

The experiment was run by controlling the arm to move from one position to another and back and so forth, following the trajectories in \figref{fig:3d-arm-desired}. Tiles of the experiment are shown in \figref{fig:robotarmvalk}.
The system identification procedure was brought online and was able to quickly identify the regressor parameters of the model.
The initial guess was purposefully chosen to be incorrect to show the convergence properties (see \figref{fig:3d-arm-r2}) of the procedure as the nominal system model -- that estimated from engineering software -- was known with reasonable accuracy.

In order to avoid bad and potentially dangerous behavior, a threshold of $.95$ was set on the coefficient of determination, $R^{2}$, below which the nominal model would be used in place of the identified model.
The validity of the identified model (and thereby its use in the controller) was also contingent upon the number of data points recorded.
Specifically, it was required that the historical data buffers be full -- in this experiment, 50 historical data were used in the buffer.
The identification procedure was performed at 3 Hz and thus it took about 16--17 seconds for the buffer to be filled.
It is quite apparent from \figref{fig:3d-arm-errors}, as noted in the caption, when the controller began to use the identified model instead of the nominal model.

%
%
%
%



\section{Conclusions}

Identification and control of a physical system, in particular an n-DOF robotic system with the implementation of an efficient computational mechanism was shown and demonstrated on three rigid body manipulators. The regressor involved in this implementation required a run time of $O(n^2)$ and computational errors resulting from this algorithm are solely due to the error in measurement of the states of the robots. 
This is a numerical method for computing the regressor and does not use symbolic expressions which are important for a robot like Valkyrie which has 44 degrees of freedom. Since SVA is required for computed torque control, evaluating the regressor through Algorithm 1 requires no extra computational overhead. 
In other words, this procedure can be directly integrated within the Rigid Body Dynamics Library \cite{featherstone:RBDL}. 
%


\bibliographystyle{plain}
\bibliography{Styles/bibdata}

\end{document}